

\catcode`\@=11
\def\undefine#1{\let#1\undefined}
\def\newsymbol#1#2#3#4#5{\let\next@\relax
 \ifnum#2=\@ne\let\next@\msafam@\else
 \ifnum#2=\tw@\let\next@\msbfam@\fi\fi
 \mathchardef#1="#3\next@#4#5}
\def\mathhexbox@#1#2#3{\relax
 \ifmmode\mathpalette{}{\m@th\mathchar"#1#2#3}%
 \else\leavevmode\hbox{$\m@th\mathchar"#1#2#3$}\fi}
\def\hexnumber@#1{\ifcase#1 0\or 1\or 2\or 3\or 4\or 5\or 6\or 7\or 8\or
 9\or A\or B\or C\or D\or E\or F\fi}

\newdimen\ex@
\ex@.2326ex
\def\varinjlim{\mathop{\vtop{\ialign{##\crcr
 \hfil\rm lim\hfil\crcr\noalign{\nointerlineskip}\rightarrowfill\crcr
 \noalign{\nointerlineskip\kern-\ex@}\crcr}}}}
\def\varprojlim{\mathop{\vtop{\ialign{##\crcr
 \hfil\rm lim\hfil\crcr\noalign{\nointerlineskip}\leftarrowfill\crcr
 \noalign{\nointerlineskip\kern-\ex@}\crcr}}}}
\def\varliminf{\mathop{\underline{\vrule height\z@ depth.2exwidth\z@
 \hbox{\rm lim}}}}

\font\tenmsa=msam10
\font\sevenmsa=msam7
\font\fivemsa=msam5
\newfam\msafam
\textfont\msafam=\tenmsa
\scriptfont\msafam=\sevenmsa
\scriptscriptfont\msafam=\fivemsa
\edef\msafam@{\hexnumber@\msafam}
\mathchardef\dabar@"0\msafam@39
\def\dashrightarrow{\mathrel{\dabar@\dabar@\mathchar"0\msafam@4B}}
\def\dashleftarrow{\mathrel{\mathchar"0\msafam@4C\dabar@\dabar@}}

\font\tenmsb=msbm10
\font\sevenmsb=msbm7
\font\fivemsb=msbm5
\newfam\msbfam
\textfont\msbfam=\tenmsb
\scriptfont\msbfam=\sevenmsb
\scriptscriptfont\msbfam=\fivemsb
\edef\msbfam@{\hexnumber@\msbfam}
\def\Bbb#1{{\fam\msbfam\relax#1}}
\def\widehat#1{\setbox\z@\hbox{$\m@th#1$}%
 \ifdim\wd\z@>\tw@ em\mathaccent"0\msbfam@5B{#1}%
 \else\mathaccent"0362{#1}\fi}
\font\teneufm=eufm10
\font\seveneufm=eufm7
\font\fiveeufm=eufm5
\newfam\eufmfam
\textfont\eufmfam=\teneufm
\scriptfont\eufmfam=\seveneufm
\scriptscriptfont\eufmfam=\fiveeufm

\newsymbol\boxtimes 1202

\catcode`\@=12

\magnification=\magstep1
\font\title = cmr10 scaled \magstep2

\font\smalltext = cmr7
\font\smallmath= cmmi7
\font\tinymath=cmmi5
\font\smallsym = cmsy7
\font\author = cmcsc10
\font\addr = cmti7
\font\byabs = cmr7

\parindent=1em
\baselineskip 15pt
\hsize=12.3 cm
\vsize=18.5 cm

\newcount\refcount
\newcount\seccount
\newcount\sscount
\newcount\eqcount
\newcount\boxcount
\newcount\testcount
\newcount\bibcount
\boxcount = 128
\seccount = -1

\def\proc#1#2{\advance\sscount by 1
	\medskip\goodbreak\noindent{\author #1}
	{\tenrm{\number\sscount}}:\ \ {\it #2}}
\def\nproc#1#2#3{\advance\sscount by 1\global
	\edef#1{#2\ \number\sscount}	
	\medskip\goodbreak\noindent{\author #2}
	{\tenrm{\number\sscount}}:\ \ {\it #3}}
\def\proof{\medskip\noindent{\it Proof:\ \ }}

\def\eql#1{\global\advance\eqcount by 1\global
	\edef#1{(\number\eqcount)}\leqno{#1}}
\def\ref#1#2{\advance\refcount by 1\global
	\edef#1{[\number\refcount]}\setbox\boxcount=
	\vbox{\item{[\number\refcount]}#2}\advance\boxcount by 1}
\def\biblio{{\frenchspacing
	\bigskip\goodbreak\centerline{\bf REFERENCES}\medskip
	\bibcount = 128\loop\ifnum\testcount < \refcount
	\goodbreak\advance\testcount by 1\box\bibcount
	\advance\bibcount by 1\vskip 4pt\repeat\medskip}}

\def\emph{\bf}
\def\colon{{:}\;}
\def\|{|\;}


\def\scirc{{\scriptstyle\circ}}
\def\rdn#1{\downarrow\rlap{$\scriptstyle #1$}}
\def\ldn#1{\llap{$\scriptstyle #1$}\downarrow}

\def\R{{\Bbb R}}
\def\Z{{\Bbb Z}}

\def\F{{\Bbb F}}

\def\A{{\Bbb A}}

\def\G{{\cal G}}

\def\qed{\hfill\hbox{$\sqcup$\llap{$\sqcap$}}\medskip}

\def\P{{\Bbb P}}
\def\SL{{\rm SL}}
\def\Sp{{\rm Sp}}
\def\SO{{\rm SO}}
\def\SU{{\rm SU}}
\def\Spin{{\rm Spin}}
\def\PSL{{\rm PSL}}
\def\PGL{{\rm PGL}}
\def\Spec{{\rm Spec}\,}

\ref\Borel{Borel, A.: On free subgroups of semisimple groups.  
	{\it Enseign. Math. (2)} {\bf 29}  (1983),  no. 1-2, 151--164.}

\ref\EGA{Dieudonn\'e, J.; Grothendieck, A.: \'El\'ements
	de G\'eom\'etrie Alg\'ebrique IV, {\it Publ. Math. I.H.E.S.} {\bf 20} (1964), 
	{\bf 24} (1965), {\bf 28} (1966), {\bf 32} (1967).}

\ref\LPink{Larsen, M.; Pink, R.: Finite subgroups of algebraic groups, preprint.}


\ref\Tits{Tits, J.: Free subgroups in linear groups, {\it J. Algebra} {\bf 20} (1972) 250--270.}

\centerline{\title WORD MAPS HAVE LARGE IMAGE}
\bigskip
\centerline{\byabs BY}
\medskip
\noindent{\author\hfill Michael Larsen\footnote*
{\tenrm Partially supported by NSF
Grant DMS-0100537}\hfill}
\medskip
\centerline{\addr Department of Mathematics, Indiana University}
\centerline{\addr Bloomington, IN 47405, USA}
\bigskip

\centerline{\byabs ABSTRACT}
\smallskip
{\byabs \narrower\narrower
\textfont0 = \smalltext
\textfont1 = \smallmath
\scriptfont1 = \tinymath
\textfont2 = \smallsym
An element $w$ in the free group on $r$ letters defines
a map $f_{w,G}\colon G^r \to G$
for each group $G$.  In this note, we show that whenever $w\neq 1$
and $G$ is a semisimple algebraic group, $f_{w,G}$ is dominant.  As an application,
we show that for fixed $w$ and $\Gamma_i$ a sequence of pairwise non-isomorphic 
finite simple groups,
$$\lim_{i\to\infty}{\log |\Gamma_i|\over\log |f_{w,\Gamma_i}(\Gamma_i^r)|} = 1.$$
}

\medskip
Let $F_r$ be the free group on $r$ generators $x_1,\ldots,x_r$.
For any group $G$, each word
$$w=x_{a_1}^{b_1}x_{a_2}^{b_2}\cdots x_{a_m}^{b_m}\in F_r$$
defines a corresponding {\emph word map} $f_{w,G}\colon G^r\to G$:
$$f_{w,G}(g_1,\ldots,g_r) = g_{a_1}^{b_1}g_{a_2}^{b_2}\cdots g_{a_m}^{b_m}.$$
The main result of this note is as follows:

\nproc\Dominant{Theorem}{If $G$ is a simple algebraic group over any field $K$ and $w\neq 1$,
then $f_{w,G}$ is a dominant morphism.  In other words, $f_{w,G}(G)$ contains a non-trivial Zariski-open subset of $G$.}
\medskip
As an application, we prove the following theorem, which answers a question of A.~Shalev:

\nproc\Plentiful{Theorem}{If $w\neq 1$ and $\Gamma_1,\Gamma_2,\ldots$ is an infinite sequence
of finite simple groups, no two isomorphic to one another, then
$$\lim_{i\to\infty}{\log |\Gamma_i|\over\log |f_{w,\Gamma_i}(\Gamma_i^r)|} = 1.$$
}

I would like to thank A.~Shalev for some interesting discussions regarding this paper
and acknowledge the hospitality of the Hebrew University, where most of the work
was carried out.

Shortly after this paper appeared in the {\it Israel Journal of Mathematics}, I learned from 
V.~Platonov that A.~Borel \Borel\ discovered \Dominant\ twenty years ago.  As far as I am aware,
the application, \Plentiful, is new.

Without loss of generality, we may assume $K$ is algebraically closed.
If $\pi\colon G\to H$ is any morphism of algebraic groups, the diagram
$$
\matrix{
G^r&{\buildrel f_{w,G}\over\longrightarrow}&G\cr
\ldn{\pi^r}&&\rdn{\pi}\cr
H^r&{\buildrel f_{w,H}\over\longrightarrow}&H\cr
}\eql\square
$$
commutes.  Applying \square\ when $\pi$ is an isogeny, we see that it suffices to prove the
theorem for $G$ simply connected.  Applying it to the factor inclusion maps
when $G$ is a product, we see that it suffices to consider the case of simply connected
almost simple groups.  Such groups are indexed by connected Dynkin diagrams, and
we begin with type A.

\nproc\An{Lemma}{\Dominant\ holds for $G=\SL_n$.}

\proof
We use induction on $n$, the base case $n=1$ being trivial.  Define $\chi_n\colon\SL_n\to\A^{n-1}$
so that if $g\in\SL_n$ has characteristic polynomial
$$x^n-a_1x^{n-1}+a_2 x^{n-2}-\cdots+(-1)^n,$$
then
$$\chi_n(g)=(a_1,a_2,\ldots,a_{n-1}).$$
Thus $\chi_n$ is constant on conjugacy classes of $\SL_n$.  Over the non-empty open
subvariety of $\A^{n-1}$ corresponding to polynomials with non-zero discriminant, the fibers
of $\chi_n$ are single conjugacy classes.  Since $f_{w,\SL_n}(\SL_n^r)$
is a union of conjugacy classes of $\SL_n$, it contains a dense open subset
of $\SL_n$ if and only if its image under $\chi_n$ contains a dense
open subset of $\A^{n-1}$.

The induction hypothesis and the inclusion $\SL_{n-1}\hookrightarrow\SL_n$
imply that the Zariski closure of the image of $\chi_n\scirc f_{w,\SL_n}$
contains a dense open subset of the hyperplane
$$\{(a_1,\ldots,a_{n-1})\mid 1-a_1+a_2-\cdots+(-1)^n=0\}\eql\plane$$
corresponding to elements of $\SL_n$ with eigenvalue $1$.  On the other hand,
$\SL_n$ is connected, so the Zariski-closure of
$$\chi_n(f_{w,\SL_n}(\SL_n^r))$$
is connected.  To prove that the closure is all of $\A^{n-1}$
(and therefore that $f_{w,\SL_n}$ is dominant), we need only show that 
some element of the image of $\chi_n\scirc f_{w,\SL_n}$
is not contained in \plane, i.e., that some element of
$\SL_n$ in the image of the word map does not have 1 as an eigenvalue.

To do this, we begin with a global field $F$ contained in $K$.  Let $D$ be a division algebra
of degree $n$ over $F$ and $\SL_1(D)$ the multiplicative group of
elements of $D^\times$ with reduced norm $1$, which we regard as the group of $F$-points
of an inner form $S$ of $\SL_n$ over $F$.  Let $x_1\in\SL_1(D)=S(F)$ denote an
element of infinite order and $x_1,\,x_2,\,x_3,\,\ldots$ a maximal sequence of
elements in $S(F)$ such that $x_{n+1}$ does not lie in the normalizer of the identity 
components of the Zariski closure $X_n$ of the subgroup generated by $x_1,\,\ldots,\,x_n$.
Such a sequence is finite since $\dim X_{n+1}>\dim X_n$.  Let $\Gamma$
be the subgroup of $\SL_1(D)$ generated by all the $x_i$.  As $\Gamma$ is finitely generated
and Zariski-dense in the semisimple group $S$,
the Tits alternative \Tits\ implies it contains a subgroup isomorphic to $F_r$.
The inclusions
$$F_r\subset\Gamma\subset\SL_1(D)\subset D$$
allow us to regard $w$ as an element of $D\setminus\{1\}$.  In particular, $w-1\in D$
is non-zero and so invertible in $D$.  As $K$ is algebraically closed,
$S(K) = \SL_n(K)$, so it follows that
$f_{w,\SL_n}(\SL_n(K)^r)$ contains an element of the desired kind.\qed

At this point we know that $f_{w,G}$ is dominant for any semisimple group $G$ whose 
Dynkin diagram components are all of type A.  Suppose $G$ is a group of this type
and $G\hookrightarrow H$ is an injective homomorphism of semisimple groups of equal rank;
that is, a maximal torus $T$ of $G$ is again a maximal torus of $H$.  Then 
the image of $f_{w,G}$ contains a dense open subset of $T$.  Let $\Psi\colon H\times H\to H$
be the conjugation morphism defined by 
$$\Psi(h_1,h_2)=h_1h_2h_1^{-1}.$$
Then,
$$f_{w,H}(H^r)\supset\Psi(H\times f_{w,G}(G^r))\supset \Psi(H\times (T\cap f_{w,G}(G^r))).$$
The restriction of $\Psi$ to $H\times T$ is dominant since every semisimple element of $H$
is conjugate to an element of $T$.  Therefore, $f_{w,H}(H^r)$ is dense in $H$. 
Thus, we need only verify:

\proc{Lemma}{Every simply connected almost simple Lie group $H$ over $K$ contains an 
equal rank semisimple subgroup whose Dynkin diagram components are all of type A.}

\proof
There are obvious inclusions
$$\SL_2^n\subset \Sp_{2n},$$
$$\SL_2^{2n}=\Spin_4^n\subset \Spin_{4n}\subset\Spin_{4n+1},$$
$$\SL_2^{2n-2}\times\SL_4=\Spin_4^{n-1}\times\Spin_6\subset\Spin_{4n+2}\subset\Spin_{4n+3}.$$
For the exceptional groups we use the fact that a closed root subsystem
of an irreducible root system gives rise to a reductive subgroup;
if the subsystem has equal rank, the same will be true on the group level.  We apply this to
the (root system) inclusions
$$A_2^3\subset E_6,\ A_1\times A_3^2\subset E_7,\ A_4^2\subset E_8,
\ A_2^2\subset F_4,\ A_2\subset G_2$$
to prove the lemma.\qed

This finishes the proof of \Dominant.\qed

\proc{Corollary}{If $G$ is a  semisimple algebraic group and $w$ is a non-trivial
element of $F_r$, then $f_{w,G(\R)}(G(\R)^r)$ has non-empty interior.}

\proof By definition 
(\EGA~IV~17.3.7), smoothness of morphisms is an open property, 
so generic smoothness of
a morphism of varieties can be checked at the generic point, where
it is equivalent to separability of the extension 
of function fields (\EGA~$0_{\rm IV}$~19.6.1).
As $f_{w,G}$ is a dominant morphism between varieties in characteristic 0,
$G^r$ contains a non-empty open subvariety of smooth points.
As $G(\R)$ is Zariski-dense in $G$, there exists a smooth point 
$x\in G^r(\R)$.  The image
of $x$ in $G(\R)$ is an interior point by the implicit function theorem.
\qed

\proc{Question}{Is $f_{w,G}$ always surjective at the algebraic variety level?  How about
at the level of $\R$-points?}
\medskip

Finally, we prove \Plentiful.
We use the classification of finite simple groups to divide the
problem into three parts:
groups of Lie type of bounded dimension, classical groups in the limit as rank tends to $\infty$,
and alternating groups.  We can disregard sporadic groups because we are interested only
in behavior in the limit.  We begin with the part of the problem directly related to the
algebraic group case.

\nproc\Bounded{Proposition}{For any non-trivial word $w$ and any root system $\Phi$, there
exists a constant $c>0$ such that for all simple groups $\Gamma$ of Lie type 
associated to the root system $\Phi$,
$$|f_{w,\Gamma}(\Gamma^r)| > c|\Gamma|.$$
}

\proof
The idea is to find an upper bound on the size of the fibers of $f_{w,\Gamma}$ by regarding
them, more or less, as the $\F_q$-points of fibers of a morphism of varieties 
$f_{w,G}\colon G^r\to G,$ where $G$ is a simple algebraic group with root system $\Phi$
and $\F_q$ is a finite field.  The basic estimate is the naive one:
$$|X(\F_q)| < C q^{\dim X},$$
but there are a number of technical difficulties in making this strategy work.
To begin with, it is not quite accurate to identify $\Gamma$ with a group of
the form $G(\F_q)$.  This is especially problematic when $\Gamma$ is a Suzuki or
Ree group.  The constant $C$ above has to be uniform across fibers of $f_{w,G}$
and independent of characteristic.  Although by \Dominant,
generically the fibers of $f_{w,G}$ have dimension $(r-1)(\dim G)$, some fibers
may have higher dimension, and we must account for these.  Rather than developing from scratch
a technology to deal with these problems, we appeal to \LPink, where such a technology
already exists.

Let $G$ be an adjoint simple group with root system $\Phi$ over an algebraically
closed field $K$ of characteristic $p$ and $\Gamma\subset G(K)$.  Without loss of generality,
$|\Gamma|\gg 0$, so by \LPink\ Prop. 3.5, $\Gamma$ is sufficiently general.  
The dimension of fibers of
$f_{w,G}$ is upper semicontinuous, so there exists a proper closed reduced subscheme
$X_{w,G}\subset G^r$ such that $f_{w,G}$ restricted to $G^r\setminus X_{w,G}$ has 
constant fiber dimension.

The subscheme $X_{w,G}$ depends on the characteristic $p$, but the set of 
all such subschemes forms
a constructible family in $\G^r$, where $\G/\Spec\Z$ is the adjoint Chevalley scheme with
root system $\Phi$.  By \LPink\ Th. 4.3,
$$|X_{w,G}\cap\Gamma^r| < c_1 |\Gamma|^{\dim X_{w,G}\over\dim G}
\le c_1|\Gamma|^{r-{1\over \dim G}}.$$
The fibers $Y_g$ of $f_{w,G}$ as $g$ ranges over $G$ and $p$ ranges over all prime numbers
again form a constructible family, so
$$|Y_g\cap\Gamma^r| < c_2|\Gamma|^{\dim Y_g\over\dim G}.$$
Therefore,
$$|f_{w,\Gamma}(\Gamma^r)|\ge {|\Gamma|^r-|X_{w,G}\cap\Gamma^r|\over c_2|\Gamma|^{r-1}}
> {1\over c_2}|\Gamma|\left(1-c_1|\Gamma|^{-{1\over\dim G}}\right)>{|\Gamma|\over 2c_2}$$
for $|\Gamma|\gg 0$.
\qed

\nproc\Alt{Proposition}{Let $A_n$ denote the alternating group on $n$ letters.  
Then for all $\epsilon>0$ there exists $N$ such that
$$|f_{w,A_n}(A_n^r)|\ge |A_n|^{1-\epsilon}$$
for all $n\ge N$.
}

\proof
Let $\phi(n)$ and $\tau(n)$ denote the Euler $\phi$-function and the number of divisors 
function respectively.  For any $\epsilon > 0$ and any sufficiently large prime power $q$,
$\tau(q)<q^\epsilon$; as $\tau(n)$ is multiplicative, $\tau(n)=o(n^\epsilon)$.
Therefore, the number of elements of order $<n^{1-\epsilon}$ in $\Z/n\Z$  is
$$\sum_{d\vert n\atop d < n^{1-\epsilon}}\phi(d)\le \tau(n) n^{1-\epsilon}
< n^{\epsilon/2}(n^{1-\epsilon}) < n^{1-\epsilon/2}$$
for $n\gg 0$.

Now, let $p$ be an odd prime, $\Gamma=\PSL_2(\F_p)$, $I_p = f_{w,\Gamma}(\Gamma^r)$.
We claim that for $p\gg 0$, $I_p$ contains an element of order $>p^{1-\epsilon}$.
Let $\Delta_1\cong \Z/{p-1\over 2}\Z$ (resp. $\Delta_2\cong\Z/{p+1\over 2}\Z$) 
denote subgroups of $\Gamma$ associated to a split 
(resp. non-split) torus in $\PGL_2$.  Two elements $x,y\in\Delta_i$ are conjugate in $\Gamma$
if and only if $x=y^{\pm1}$.  The centralizer of any non-identity element of $\Delta_i$
is $\Delta_i$ itself, so the conjugacy class of such an element has order $|\Gamma|/|\Delta_i|$,
and no such conjugacy class meets both $\Delta_1$ and $\Delta_2$.  We conclude that 
if $p\gg 0$, the
set of elements in $\Gamma$ conjugate to some element of order $< p^{1-\epsilon}$ is at most
$|\Gamma| p^{-\epsilon/2}$.  By \Bounded, $I_p$ contains an element of order $\ge p^{1-\epsilon}$
if $p\gg 0$.

Next we consider the action of $\Gamma$ on the finite projective line $\P^1(\F_p)$
(by fractional linear transformations).  This gives an embedding
$$\PSL_2(\F_p)\hookrightarrow A_{p+1}.$$
A non-identity element of $\Delta_1$ (resp. $\Delta_2$) fixes 2 (resp. 0)
points of $\P^1(\F_p)$; if its order is $d$, its image in $A_{p+1}$ consists of
$p-1\over d$ (resp $p+1\over d$) $d$-cycles and 2 (resp. 0) 1-cycles.  Let
$$S=\{p+1\mid p\hbox{ prime}\}.$$
By the prime number theorem, the greedy algorithm guarantees that there exists an integer $B$
such that every interval of length $B$ in the set of positive integers contains the sum of
a sequence of elements of $S$, each larger than the sum of all that come after.  
In other words, for every positive integer $n$,
$$A_n\supset A_{p_1+1}\times\cdots\times A_{p_k+1}\supset \PSL_2(\F_{p_1})\times\cdots\times
\PSL_2(\F_{p_k}),$$
where
$$n-B\le p_1+1+\cdots+p_k+1\le n,\ p_i+1\in S,\ k\le\log_2 n.$$
It follows that $f_{w,A_n}(A_n^r)$ has an element which decomposes in $c = O(n^\epsilon\log n)$
cycles (including cycles of length $1$).  
The centralizer of a product of $c$ cycles in $S_n$ has order $\le n^c c! = o( |A_n|^\epsilon)$
for $n\gg 0$.  Therefore, $f_{w,A_n}(A_n^r)$ contains a conjugacy class with more than
$|A_n|^{1-\epsilon}$ elements.
\qed

\proc{Proposition}{For all $w\neq 1$ and $\epsilon > 0$ there exists $N$ such that if
$\Gamma$ is a finite simple group of Lie type of rank $>N$ then
$$|f_{w,\Gamma}(\Gamma^r)|>|\Gamma|^{1-\epsilon}.$$
}

\proof 
Suppose $\Gamma$ has a central extension $\tilde\Gamma$ in the set of groups
$$\{\SL_n(\F_q),\SO_{n,n}(\F_q),\SO_{2n+1}(\F_q),\SO_{n+2,n},\Sp_{2n}(\F_q),
\SU_{2n}(\F_q),\SU_{2n+1}(\F_q)\}.\eql\list$$
We note the inclusions
$$A_n\subset\SL_n(\F_q)\subset\SO_{n,n}(\F_q)\subset\SO_{2n+1}(\F_q)\subset\SO_{n+2,n}(\F_q)
\subset\SL_{2n+2}(\F_q),\eql\Ortho$$
$$A_n\subset\SL_n(\F_q)\subset\Sp_{2n}(\F_q)\subset\SL_{2n}(\F_q),\eql\Symp$$
and
$$A_n\subset\SL_n(\F_q)\subset\SU_{2n}(\F_q)\subset\SU_{2n+1}(\F_q)\subset\SL_{4n+2}(\F_q).\eql\Unit$$
Suppose there exists a constant $\delta < 1$ (depending on $w$) such that
for all $n\gg 0$ there exists $x\in f_{w,A_n}(A_n^r)\subset A_n$
whose images in $\SL_m(\F_q)$ for 
$$m\in \{2n+2, 2n, 4n+2\}$$
all have centralizer orders
$O(q^{n^{1+\delta}})$.   This implies the same upper bound for the centralizer
of the image $y$ of $x$ in $\tilde\Gamma$.
The order of $\tilde\Gamma$ is at least
$$|\SL_n(\F_q)|= {1\over q-1}\prod_{i=0}^{n-1} (q^n-q^i)>q^{n^2-1}\prod_{j=2}^\infty(1-q^{-j})
>{q^{n^2-1}\over 2},$$
so the conjugacy class of $y$ in $\tilde\Gamma$ has order at least
$|\tilde\Gamma|^{1-\epsilon}$ if $n\gg 0$.  
In mapping from $\tilde\Gamma$ to $\Gamma$
the size of a conjugacy class goes down by at most a factor of 2n+1.
The estimate $O(q^{n^{1+\delta}})$ is therefore enough to prove the proposition.

The composed maps $A_n\subset\SL_m(\F_q)$ in \Ortho, \Symp, and \Unit\ factor through 
$A_{2n+2}$, $A_{2n}$, and $A_{4n+2}$ respectively, and an element in $A_n$
consisting of $c$ cycles maps to an element with $2c+2$, $2c$, and $4c+2$ cycles respectively.
As in \Alt, we can find $z\in A_m$, the image of $x\in A_n$, such 
that $z$ consists of $O(m^\epsilon\log m)$ cycles.  Regarding $z$ as a permutation matrix
in $\SL_m(\F_q)$, we consider its centralizer in the matrix algebra $M_m(\F_q)$.
If $(a_{i,j})$ is a matrix commuting with the permutation matrix associated with a permutation
$\sigma$, then
$$a_{i,j}=a_{\sigma(i),\sigma(j)}$$
for all $i,j$.  Therefore, the whole matrix is determined by any set of rows representing
all $\sigma$-orbits.  If $\sigma$ has $\le 4c+2$ orbits, the centralizer has order $\le q^{(4c+2)m}$.  
Therefore the centralizer of the image of $z$ in $\SL_m(\F_q)$ (or in any subgroup thereof) 
has order $O(q^{m^{1+2\epsilon}})$.  The proposition, and therefore \Plentiful, follows.
\qed

\biblio
\end